\documentstyle[12pt]{amsart}

\newtheorem{theorem}{Theorem}[section]
\newtheorem{lemma}[theorem]{Lemma}
\newtheorem{corollary}[theorem]{Corollary}

\theoremstyle{definition}

\newtheorem{definition}[theorem]{Definition}
\newtheorem{conjecture}[theorem]{Conjecture}
\newtheorem{example}[theorem]{Example}

\newtheorem{remark}[theorem]{Remark}

\def\Cal{\cal}

\def\fix{f}

\begin{document}
\title{Canonical heights and entropy
in arithmetic dynamics}
\author{M. Einsiedler}
\address{(M.E.) Mathematical Institute, University of Vienna,
Strudlhofgasse 4, A-1090 Wien, Austria.}
\email{manfred@@mat.univie.ac.at}
\author{G. Everest}
\author{T. Ward}
\address{(G.E. \& T.W.) School of Mathematics, University of East Anglia,
Norwich NR4 7TJ, UK.}
\email{g.everest@@uea.ac.uk}
\dedicatory{\today}
\thanks{The first author acknowledges the support of
EPSRC postdoctoral award GR/M49588, the
second thanks Jonathan Lubin and
Joe Silverman for the AMS Sectional meeting on Arithmetic
Dynamics at Providence, RI, 1999}
\begin{abstract}
The height of an algebraic number in the sense of
Diophantine geometry is known to be related to
the entropy of an automorphism of a solenoid associated to the
number. An elliptic analogue
is considered,
which necessitates introducing a notion of entropy for
sequences of transformations.
A sequence of transformations are defined
for which there is a canonical arithmetically defined
quotient whose entropy is the canonical height, and
in which the fibre entropy is accounted for by
local heights at primes of singular reduction, yielding
a dynamical interpretation of singular reduction.
This system is related to local systems, whose entropy
coincides with the local canonical height up to sign.
The proofs use transcendence theory, a strong form
of Siegel's theorem, and an elliptic analogue
of Jensen's formula.
\par
These elliptic systems are based upon iteration
of the duplication map; the ideas extend
to morphisms of projective space, giving examples
where the associated entropies coincide with the morphic 
heights of Call and Goldstine.
In particular, the local morphic heights
at infinity for polynomials are realized as integrals over an associated
Julia set with respect
to the maximal measure, giving an analogue of the
Jensen formula in that setting also.
\end{abstract}
\maketitle
\section{Introduction}
Let $Q$ denote a finite rational point of the projective line
${\Bbb P}^1$. Then $Q$ has an associated dynamical system 
$T_Q:X_Q\rightarrow X_Q$, where $T_Q$ is a continuous map on
an underlying compact group $X_Q$
known as a
\it{solenoid }\rm (defined later).  
The topological
entropy of this dynamical
system, an intrinsic invariant measuring orbit
complexity, coincides with the
Diophantine height $h(Q)$ of $Q$.
If $Q=[q,1]$ corresponds to the rational number
$q=a/b$ then
this height can be written using Jensen's formula
as an integral,
\begin{equation}\label{jensen}
h(Q)=\log\max\{\vert a\vert,\vert b\vert\}=\int_{\Bbb T}\log |bx-a|dm,
\end{equation}
where $\Bbb T$ is the unit circle and $m$ is Haar measure.
The number of elements of $X_Q$ fixed by $T_Q^n$
is $|b^n-a^n|$. Writing $\phi_n(x)=x^n-1$, the
polynomial whose roots form the $n$-torsion subgroup of
the unit circle, gives $|b^n-a^n|=|b^n\phi_n(a/b)|$.
The main point of reference is the approach taken
in \cite{lind-ward-1988}, where the entropy
is calculated by noting that the space $X_Q$ is covered by the adeles
and the dynamics lift nicely. The lifted map restricts
to the local components, and the local entropies agree
with the local projective heights. In this covering space,
the periodic point data is destroyed however.

The arithmetic side of the last paragraph has a direct
analogue in which (roughly speaking) ${\Bbb T}$ is
replaced by a complex elliptic curve, and the projective
height is replaced by the global canonical height, which
is known to decompose as a sum of local canonical heights. 
The cyclotomic division polynomials carrying knowledge
of torsion in the circle are replaced by
the elliptic division polynomials.
Several attempts
have been made to find the fourth corner of this
square of ideas, namely a family of {\em elliptic dynamical
systems}, whose topological entropy is given by the
canonical height on the curve, and whose periodic point data
is given by expressions involving the elliptic division
polynomial (see \cite{dambros-everest-miles-ward},
\cite{everest-ward-1998}, \cite{everest-ward-1999}).

In this paper, we have three objectives.
The first is to show that by widening the usual concept
of a dynamical system to include
sequences of transformations which are not necessarily
the iterates of a single transformation, we can construct
dynamical systems from rational points on elliptic curves
which interpret the known arithmetic properties
of heights. The results include a dynamical
interpretation of the phenomenon of singular reduction.
In this wider concept of dynamics, there is
a natural notion of entropy, measuring growth in orbit
complexity along the sequence. The maps act on the adeles,
just as in \cite{lind-ward-1988}; they are 
built
from the duplication map on the underlying elliptic curve.
Duplication can be viewed as a morphism on the projective 
line, which
informs our second objective: to construct sequences of 
transformations built from iterates of morphisms on projective
space, and relate the entropy to the canonical height for these
morphisms. The underlying space for all of these systems
is locally compact (it is the adele ring). 
The third objective is to argue through examples that
if the underlying
space is compact this places severe restrictions upon
the volume growth rates in our entropy calculations. 
By this we mean to imply that our maps on the adeles
are natural from the point of view of the systems we seek.

Since we are bringing together two areas (arithmetic and dynamics)
with a fair amount of technical detail, the
main conclusions are stated with precise definitions later. 
In
\cite{lind-ward-1988}, the entropy was calculated by
showing it is equal to that of the diagonal
multiplication by $q$ on the
rational adeles.

\vskip.1in
\noindent{\bf Theorem} (see Section \ref{dynintofbadred})
{\it Let $E$ denote an
elliptic curve defined over $\Bbb Q$, and
$Q$ a rational point on $E$. Then
$Q$ generates a sequence of diagonal transformations
$\bold U$ on the adeles with
the following properties.

\noindent1. If $Q$ has non-singular reduction modulo $p$ for all primes $p$
then the entropy $h({\bold U})=\hat h(Q)$, the global canonical height of $Q$.

\noindent2. Let $S$ denote the set of primes $p$ for which $Q$ has singular
reduction modulo $p$; write 
${\Bbb Q}_{S}=\prod_{p\in S}{\Bbb Q}_p$,
and ${\bold U}_{S}$ for the restriction of
$\bold U$ to 
${\Bbb Q}_{S}$. Then the quotient entropy
$h({\bold U}/{\bold U}_S)=\hat h(Q)$.}
\vskip.1in

Duplication on an elliptic curve provides a natural
example of a morphism on ${\Bbb P}^1$ of degree 4. Suppose now that 
$F$ denotes an
arbitrary morphism
on ${\Bbb P}^1$ defined over $\Bbb Q$,
corresponding to the rational function $f$ in
one variable.
In \cite{call-goldstine-1997} a notion of canonical height
is attached to $F$.
Write $\hat h_f(Q)$ for the global canonical height of 
$Q\in {\Bbb P}^1(\Bbb Q)$ and, for each
prime $p\le\infty$, write $\lambda_{f,p}(Q)$ for the local canonical height.
If $f$ is a polynomial, we
obtain a simple result which matches up local and global
heights with the local and global entropies attached to a natural
sequence of transformations coming from the iterates of $f$ on
a single rational number. This point of view is closer
to dynamical systems in the usual sense: iterates
of rational maps provide the raw material for the
transformations.

\vskip.1in
\noindent{\bf Theorem}
{(see Section \ref{goldstine})}
{\it Suppose $F:{\Bbb P}^1({\Bbb Q})\to
{\Bbb P}^1({\Bbb Q})$ corresponds
to a polynomial $f$ in one variable. Let $Q\in {\Bbb P}^1(\Bbb Q)$ 
denote a finite rational point corresponding to $q\in \Bbb Q$.
Then the iterates of $f$ on $q$ generate a sequence of 
transformations ${\bold T}$
on the adeles with the following properties:

\noindent1. The global entropy $h({\bold T})=\hat h_f(q)$,

\noindent2. The restriction
${\bold T}_p$ of ${\bold T}$ to ${\Bbb Q}_p$ has entropy
$h({\bold T}_p)=\lambda_{f,p}(q)$.}
\vskip.1in

This system emulates an automorphisms of the
solenoid. It is natural to ask what can be done for
a general morphism -- one that does not correspond to a polynomial.
The elliptic system
provides an answer for a particular morphism
of degree 4, and clarifies some problems with the general case.
Examples show that
these more general results may be interpreted in terms of known objects.
For example, the local heights in the elliptic case, and in the
polynomial case of a morphism, arise as integrals over an
associated Julia
set. In the elliptic case, this is because the Julia set coincides
with the elliptic curve. In the case of the second result, it
is possible to see this from what is known already. For the more
general case, it is an open problem that 
completes
this circle of ideas in a satisfactory way.

\section{Definitions and background on entropy}
Most of the definitions and results below
are straightforward modifications of well-known
theory, so the
results are simply stated. The interest is in the later
examples.
Let $X$ be a `space': a
standard probability space $(X,{\Cal B},\mu)$,
a compact metric space $(X,\rho)$, or a
locally compact metric space $(X,d)$.
A {\sl sequential action} on $X$ is a
sequence ${\bold T}=\left(T_n\right)_{n\ge 1}$
of maps $T_n:X\to X$ with the property that
each $T_n$ is a $\mu$-preserving $\Cal B$-measurable map,
a continuous map, or a uniformly continuous map respectively.
One of the essential features of the elliptic
phenomena we are trying to capture is that the volume
grows at some natural rate. Let
$r:{\Bbb N}\to{\Bbb R}$ be non-decreasing
with $r(n)\nearrow\infty$.
A finite partition
$\xi$ of $(X,{\Cal B},\mu)$ is a collection
$\{A_1,\dots,A_k\}$ of $\Cal B$-measurable sets
with $\mu(\bigcup_{i=1}^kA_i)=1$ and
$\mu(A_i\cap A_j)=0$ for all $i\neq j$. The
{\sl entropy} of such a partition is
$H(\xi)=-\sum_{i=1}^{k}\mu(A_i)\log\mu(A_i)$
(with the convention that $0\log 0=0$),
and the {\sl join} of $\xi$ with another finite
partition $\eta=\{B_1,\dots,B_{\ell}\}$
is the partition $\xi\vee\eta=
\{A_i\cap B_j\mid
1\le i\le k, 1\le j\le\ell\}$.
If $T:X\to X$ is a measurable map, then
$T^{-1}\xi$ denotes the partition
$\{T^{-1}A_1,\dots,T^{-1}A_k\}$.
\begin{definition}\label{mtseqentdef}
The (measure-theoretic) sequential entropy
of $\bold T$ on $(X,{\Cal B},\mu)$ is given by
\begin{equation*}
h^{r}_{\mu}({\bold T})=
\sup_{\xi}\limsup_{n\to\infty}
\frac{1}{r(n)}
H\left(\bigvee_{j=1}^{n}T_j^{-1}\xi\right),
\end{equation*}
where the supremum is taken over all finite
partitions.
\end{definition}
\begin{example}\label{metricexamplesone}
1. Let $r(n)=n$, and let $T_j=T^j$
for all $j\ge 1$ where $T$ is a single
measure-preserving transformation. Then
$h^{r}_{\mu}({\bold T})=h_{\mu}(T)$,
the usual measure-theoretic entropy of $T$.

\noindent2. Let $r(n)=n$ again, and
let $T_j=T^{a_j}$ for a fixed increasing
sequence $A=(a_1,a_2,\dots)$. Then
$h^{r}_{\mu}({\bold T})=h_A(T)$ the
`$A$-entropy' or sequence-entropy introduced by
Kushnirenko \cite{kushnirenko-1967}
as an invariant of measure-preserving
transformations not reducible to entropy or
spectral invariants unless $T$ has positive entropy
(see \cite{krug-newton-1972}).
\end{example}

It follows from the second example that $h^{r}_{\mu}$
cannot be more functorial than the $A$-entropy $h_A$.
In particular,
the relation $h^{r}_{\mu}({\bold T}\times
{\bold S})=h^{r}_{\mu}({\bold T})+h^{r}_{\mu}({\bold S})$
does not always hold (by
\cite[Example 7]{lemanczyk-1985}),
and writing ${\bold T}^k$ for the sequence
$\left(T_j^k\right)_{j\ge 1}$,
the relation $h^{r}_{\mu}({\bold T}^k)=kh^{r}_{\mu}({\bold T})$
does not always hold (by
\cite[Example 1]{lemanczyk-1985}).

Following Bowen, we next define a topological
entropy and a volume-growth entropy for the topological
context.
Let $X$ be a compact metric space
$(X,\rho)$, write $N(U)$ for the least
cardinality of a finite subcover of an open
cover $U$, and use $\vee$ to denote the
common refinement of two open covers.
\begin{definition}\label{akmdef} The (topological) entropy
of ${\bold T}$ on $(X,\rho)$ is
\begin{equation*}
h^{r}_{top}({\bold T})=\sup_{U}\limsup_{n\to\infty}
\frac{1}{r(n)}\log N\left(\bigvee_{j=1}^{n}T_j^{-1}U
\right),
\end{equation*}
where the supremum is taken over all open covers
$U$ of $X$.
\end{definition}
\begin{example}
1. Let $r(n)=n$, and let $T_j=T^j$
for all $j\ge 1$ where $T$ is a single
continuous map on $(X,\rho)$. Then
$h_{top}^{r}({\bold T})=h_{top}(T)$, the topological
entropy of $T$ introduced in
\cite{adler-konheim-mcandrew-1965}.

\noindent2. Let $r(n)=n$, and let $T_j=T^{a_j}$
for a fixed increasing sequence $A=(a_1,a_2,\dots)$.
Then $h_{top}^{r}({\bold T})=h_{top}^A(T)$ is
the topological sequence entropy (see
\cite{dekking-1980}).

\noindent3. The directional entropy introduced by Milnor
coincides with the entropy in this sense, with
$r(n)=n$, for the sequence of transformations seen
in a strip along the chosen direction
(see \cite{milnor-camaps}).
\end{example}

As in the measure-theoretic case, it follows that
$h_{top}^{r}$ cannot be better-behaved than
the topological sequence entropy. In particular,
the relation $h^{r}_{top}({\bold T}\times
{\bold S})=h^{r}_{top}({\bold T})+h^{r}_{top}({\bold S})$
does not always hold (by
\cite[Example 5]{lemanczyk-1985}),
writing ${\bold T}^k$ for the sequence
$\left(T_j^k\right)_{j\ge 1}$,
the relation $h^{r}_{top}({\bold T}^k)=kh^{r}_{top}({\bold T})$
does not always hold (by
\cite[Example 2]{lemanczyk-1985}), and
the variational principle
$h_{top}^{r}({\bold T})=\sup_{\mu}h^{r}_{\mu}({\bold T})$
(where the supremum is taken over all probabilities
$\mu$ invariant under all the $T_j$'s) does not always hold
(see \cite[Example 3]{dekking-1980}).

Definition \ref{akmdef} is less than easy to work with, and
the calculation of topological entropy is
facilitated by Bowen's introduction of
spanning and separated sets, homogeneous measures,
and volume growth. Let now $X$ be a locally compact
metric space $(X,d)$, and assume that each $T_j$ is
uniformly continuous.
\begin{definition}\label{bowendef} Let $K\subset X$ be compact.
A set $E\subset K$ is
$(n,\epsilon)${\sl-separated} under $\bold T$ if
for any distinct points $x,y$ in $E$,
there is a $j$, $1\le j\le n$,
for which $d(T_jx,T_jy)>\epsilon$.
A set $F\subset X$ $(n,\epsilon)${\sl-spans} $K$
if, for every $x\in K$ there is a $y\in F$
for which $d(T_jx,T_jy)\le\epsilon$ for $1\le j\le n$.
Let $r_n(\epsilon,K)$ (resp. $s_n(\epsilon,K)$) denote the
largest (smallest) cardinality of a separating
(spanning) set for $K$ under $\bold T$.
Then define
\begin{eqnarray*}
h^{r}_{Bowen}({\bold T})&=&
\sup_{K}\lim_{\epsilon\searrow0}\limsup_{n\to\infty}
\frac{1}{r(n)}\log r_n(\epsilon,K)\\
&=&
\sup_{K}\lim_{\epsilon\searrow0}\limsup_{n\to\infty}
\frac{1}{r(n)}\log s_n(\epsilon,K),
\end{eqnarray*}
where the supremum is taken over all compact sets
$K\subset X$, and the coincidence of the two limits
is shown as in \cite[Lemma 1]{bowen-1971}.
\end{definition}
As in the usual case, it may be shown
that $h_{Bowen}^{r}({\bold T})=h_{top}^{r}({\bold T})$
(see \cite{bowen-periodic-1971},
\cite[Sect. 7.2]{walters}) when $(X,d)$ is
compact, and that $h_{Bowen}^{r}({\bold T})$
depends only on the uniform equivalence class of
the metric $d$ (see \cite[Proposition 3]{bowen-1971}).
\begin{definition}\label{volumedef}
Assume that each $T_j$ is a uniformly continuous
map on the locally compact metric space $(X,d)$;
write
\begin{equation*}
D_n(x,\epsilon,{\bold T})=
\bigcap_{k=1}^{n}T_k^{-1}B_{\epsilon}(T_kx)
\end{equation*}
with $B_{\epsilon}$ a metric open ball of radius $\epsilon$.
Just as in \cite[Definition 6]{bowen-1971}, call
a Borel measure $\mu$ on $X$ {\sl homogeneous}
for ${\bold T}$ if $\mu$ is finite on compact sets,
positive on some compact set, and, for every $\epsilon>0$
there exist a $\delta>0$ and a $C>0$ such that
$\mu\bigl(D_n(y,\delta,{\bold T})\bigr)\le
C\mu\bigl(D_n(x,\epsilon,{\bold T})\bigr)$
for all $n\ge 1$ and $x,y\in X$. For such a
measure, the volume-growth entropy is defined to
be
\begin{equation*}
\lim_{\epsilon\searrow0}\limsup_{n\to\infty}
-\frac{1}{r(n)}\log\mu\bigl(
D_n(x,\epsilon,{\bold T})\bigr),
\end{equation*}
which is independent of $x$ by homogeneity, and
(see \cite{bowen-1971}, {\it mutatis mutandis})
it coincides with $h_{Bowen}^{r}({\bold T})$.
\end{definition}

\begin{example}\label{simpleadelicexamples} To see Definition 
\ref{volumedef}
in an arithmetic setting, let $X$ be
the locally compact ring $\Bbb Q_{\Bbb A}$ (see
\cite[Chap. IV]{weil-1974-number} for details on
the adele ring). Write elements of the adele ring
as ${\bold x}=(x_{\infty},x_2,x_3,\dots)$, then define
(for $\alpha\in{\Bbb Q}$)
$ \alpha{\bold x}=(\alpha x_{\infty},\alpha x_2,\alpha x_3,\dots).$
Let $\mu_p$ be the Haar measure on ${\Bbb Q}_p$ ($p\le\infty$)
normalized to have $\mu_p({\Bbb Z}_p)=1$ ($p<\infty$)
and $\mu_{\infty}([0,1))=1$, and write
$\mu=\prod_{p\le\infty}\mu_p$.
It is enough to consider the neighbourhood
$B=(-1,1)\times\prod_{p<\infty}{\Bbb Z}_p$ in the following examples,
since any $\epsilon$-ball around the identity contains
the image of $B$ under an automorphism of
${\Bbb Q}_{\Bbb A}$.

\noindent1. Let $p_1,p_2,p_3,\dots$ be the rational primes in their
usual order, let $T_j({\bold x})=p_1\dots p_{j}{\bold x}$,
and let $r(n)=n\log n$.
Then it is clear that
\begin{equation}\label{todancebeneaththesettingsun}
\mu\left(
\bigcap_{j=1}^{n}T_j^{-1}B\right)=\frac{1}{p_1\dots p_{n}},
\end{equation}
so $h_{Bowen}^{r}({\bold T})=1$
(this follows from the estimate
$n\log n\ll p_n\ll n\log n$ in \cite[Theorem 4.7]{apostol-ant}).

\noindent2. Let $r(n)=n\log n$, $T_j({\bold x})=(1/p_1\dots p_{j}){\bold x}$,
and let $r(n)=n$.
Then (\ref{todancebeneaththesettingsun}) holds again,
so $h_{Bowen}^{r}({\bold T})=1$ as before.
However, in this example each `local'
entropy contribution
\begin{equation*}
\limsup_{n\to\infty}
-\frac{1}{r(n)}\log\mu_p\left(
\bigcap_{j=1}^{n}T_j^{-1}A_p\right),
\end{equation*}
where $A_p={\Bbb Z}_p$ for $p<\infty$ and $A_{\infty}=(-1,1)$,
is zero. This should be contrasted with
the usual setting, where the local entropies
sum to the global entropy (see \cite{lind-ward-1988}).

\noindent3. Let $T_j({\bold x})=\prod_{p\le j}p{\bold x}$,
where the product is over all primes less than or equal to $j$
and $r(n)=n$.
As before,
\begin{equation*}
\mu\left(
\bigcap_{j=1}^{n}T_j^{-1}B\right)=\frac{1}{\prod_{p\le n}p},
\end{equation*}
so $h_{Bowen}^{r}({\bold T})$
is positive and no larger than $2\log 2$ (see
\cite[Theorem 414]{hardy-and-wright}).

\noindent4. Let
$T_j({\bold x})=j{\bold x}$, and $r(n)=\log n$.
Then it is easy to see that
$h_{Bowen}^{r}({\bold T})=1$.

\noindent5. Let $T_j({\bold x})=j!{\bold x}$ and $r(n)=n\log n$; then in
a similar way one sees that
$h_{Bowen}^{r}({\bold T})=1$ by Stirling's formula.
\end{example}

It is well-known that `there is only one entropy',
which in our context means that $r(n)$ should always
be sub-linear or
linear in $n$ in order to allow $h^{r}$ to be
positive.
However, we have already seen in Example \ref{simpleadelicexamples}
that other rates are possible: the point is that in those examples
the underlying space is not compact.
\begin{lemma}\label{noneoncompacts}
If $X$ is compact (or of finite measure), and
$\frac{r(n)}{n}\to\infty$ as $n\to\infty$, then,
for any sequence of transformations $\bold T$,
$$h_{top}^{r}({\bold T})=h_{Bowen}^r({\bold T})=
h_{\mu}^{r}({\bold T})=0.$$
\end{lemma}
\begin{pf} Consider the measure-theoretic case first:
let $\xi$ be a finite partition. Then
$H\left(\bigvee_{j=1}^{n}T_j^{-1}\xi\right)
\le nH(\xi)$ since by \cite[Theorem 4.3]{walters}
$H(\xi\vee\eta)\le H(\xi)+H(\eta)$ for any
finite partitions $\xi$ and $\eta$, and
$H(T_j^{-1}\xi)=H(\xi)$ for all $j$.
Turning to the topological case, it is convenient
to work with $h_{top}^r$. Fix an open cover
$U$ of $X$; then
$$
N\left(\bigvee_{j=1}^{n}T_j^{-1}U\right)\le\prod_{j=1}^{n}
N(T_j^{-1}U)=N(U)^n,
$$
which gives the result.
\end{pf}
\section{Solenoids}
Suppose first that $Q=[q,1]\in{\Bbb P}^1(\bar{{\Bbb Q}})$
is a finite point on the algebraic projective line.
Then the map $x\mapsto qx$ on ${\Bbb Z}[q]$ or ${\Bbb Z}[q^{\pm1}]$ 
determines
a dual map $T_Q:X_Q\to X_Q$ on the compact abelian
dual group. Identify $X_Q$ with the dual of a subgroup of
${\Bbb Q}^d$ for
some $d$; then $T_Q$ becomes the map dual to
a rational $d\times d$ matrix $A$.
The topological entropy of $T_Q$ is given by Yuzvinskii's formula,
$h_{top}(T_Q)=\log\vert s\vert+\sum_{i}\log^{+}\vert\lambda_i\vert$,
where $s$ is the g.c.d. of the denominators of the coefficients
of the characteristic polynomial of $A$, and $\{\lambda_i\}$ are
the eigenvalues of $A$ counted with multiplicity
(see \cite{yuzvinskii-1967} for
the original derivation of this result). A more
suggestive `local-to-global' formulation of this result is
given in \cite{lind-ward-1988}: $h_{top}(T_Q)=\sum_{p\le\infty}
\sum_{i}\log^{+}\vert\lambda_{i,p}\vert_{p}$, where the inner
sum is taken over the eigenvalues of $A$ in the algebraic closure
of ${\Bbb Q}_p$, and $\vert\cdot\vert_{p}$ denotes the usual
extension of the $p$-adic valuation.

\begin{example}\label{solenoidsystemexamples}
In each case the solenoid $X_Q$ is described,
and periodic points -- points whose orbit
under the map $T_Q$ is finite -- are also discussed.

\noindent1. If $q\notin\{-1,0,1\}$ is integral, then ${\Bbb Z}[q]={\Bbb Z}$,
so the dual group $X_Q$ is the circle $\Bbb T$. The map
$T_Q$ is $x\mapsto qx$ mod $1$, and it is easy to see that
$h_{top}(T)=\log\vert q\vert$. Writing $\fix_n(T_Q)=
\{x\mid T_Q^n(x)=x\}$ for the set of points of period $n$ under
$T_Q$, we see that
$\vert\fix_n(T_Q)\vert=\vert q^n-1\vert=\vert\phi_n(q)\vert$,
where $\phi_n(x)=x^n-1$ is the
$n$th division polynomial on the circle.
Notice that $(1/n)\log\vert\fix_n(T_Q)\vert\to h_{top}(T_Q)$.

\noindent2. If $q$ is an algebraic integer (non unit-root)
of degree $d$ whose
minimal polynomial has constant coefficient $\pm 1$,
then $X_Q$ is the $d$-torus ${\Bbb T}^d$, and we may take
for $A$ the companion matrix to the minimal polynomial
of $q$. A similar argument shows that
$h_{top}(T_Q)=\sum_{i}\log^{+}\vert\lambda_i\vert$, and
$\vert\fix_n(T)\vert=\prod_{i}\vert\lambda_i^n-1\vert$.
It is still the case that
$(1/n)\log\vert\fix_n(T_Q)\vert\to h_{top}(T_Q)$,
but this is non-trivial because of the possibility
of eigenvalues with unit modulus (see \cite{everest-ward-1999}
for a detailed discussion).

\noindent3. If $q=a/b$ is a rational in lowest terms, and $X_Q$ is dual
to the group ${\Bbb Z}[q^{\pm1}]={\Bbb Z}[\frac{1}{ab}]$,
then $h_{top}(T_Q)=\sum_{p\le\infty}\log^{+}\vert\frac{a}{b}\vert_p=
\log\max\{\vert a\vert,\vert b\vert\}$ is the usual projective height
of the point $[q,1]$. Here $\fix_n(T_Q)=\vert a^n-b^n\vert=
\vert b^n\phi_n(a/b)\vert$, and again
$(1/n)\log\vert\fix_n(T_Q)\vert\to h_{top}(T_Q)$.
\end{example}

Notice that the topological entropy in each case is
given by an integral over the circle by Jensen's formula.
Yuzvinskii's formula is proved in \cite{lind-ward-1988}
using an adelic covering space: Example 
\ref{solenoidsystemexamples}.3 is a natural quotient
of the map $x\mapsto qx$ on ${\Bbb R}\times
\prod_{p\vert ab}{\Bbb Q}_p$, and the entropy may be calculated
in the covering space using the following results.
Firstly, the topological entropy of the action of $A\in M_d({\Bbb Q})$
on ${\Bbb Q}_p^d$ is given by
$h_{Bowen}(A)=\sum_{i}\log^{+}\vert\lambda_{i,p}\vert_{p}$, where the
sum is taken over the eigenvalues of $A$ in the algebraic closure
of ${\Bbb Q}_p$.
Secondly, the covering map has the same topological entropy
as the quotient map:
$h_{Bowen}^r({\Bbb Q}_p\overset{\times q}\to{\Bbb Q}_p)=
h_{top}^r(X_Q\overset{\times q}\to X_Q)$.

We therefore pursue an elliptic analogue of Yuzvinskii's
formula by considering actions on the adele ring. Notice
that Lemma \ref{noneoncompacts} shows that the quadratic
growth rates found on elliptic curves preclude the possibility
of a single homeomorphism of a compact metric space realizing
`elliptic dynamics' in a non-trivial way.

This section was deliberately formulated to reveal an underlying
genus 0: passing to genus 1 brings us to elliptic curves.

\section{Background on heights and elliptic curves}
\label{heightssection}

Let $E$ be an elliptic curve defined over the rationals,
given by a generalized Weierstrass equation
\begin{equation}
\label{genweierstrassform}
y^2+c_1xy+c_3y=x^3+c_2x^2+c_4x+c_6, 
\end{equation}
where $c_1,\dots,c_6\in{\Bbb Z}$.
For each rational prime $p$, there is a continuous
function $\lambda_p:E({\Bbb Q}_p)\to{\Bbb R}$
which satisfies the {\sl parallelogram law}
\begin{equation}\label{paralawlocal}
\lambda_p(Q+P)+\lambda_p(Q-P)=\linebreak
2\lambda_p(Q)+2\lambda_p(P)-\log
\vert x(Q)-x(P)\vert_p.
\end{equation}
If it is required that the expression $\lambda_p(Q)-
\frac{1}{2}\log\vert x(Q)\vert_p$ be bounded
as $Q\longrightarrow 0$ (the identity of $E$), then there is only one
such map, called the {\sl local canonical height}.
Note that in \cite{silverman-1994},
local heights are normalized to make them invariant
under isomorphisms: this involves adding a constant which
depends on the discriminant of $E$, the local heights
in \cite{silverman-1994} then satisfy a different form
of the parallelogram law.
For a discussion of local heights in the form used
here, see \cite{silverman-1988}.
On $E({\Bbb Q})$ the global height $\hat{h}$ can be written
as a sum of local heights -- see (\ref{sumoflocalsinellipticcaseheight})
below -- which is remarkable since there is
a more direct definition using limits of projective heights.
If $0\neq Q=[x(Q),y(Q)]\in E({\Bbb Q})$ has $x(Q)=\frac{a}{b}$,
define $h_E(Q)$ to be $\frac{1}{2}\log\max\{\vert a\vert,
\vert b\vert\}$. Then ${h}_E(Q)$ coincides with
$\frac{1}{2}h([x(Q),1])$ in the usual sense of Diophantine
geometry.
Taking the logarithmic height of the identity to
be zero gives the alternative definition
\begin{equation*}
\hat{h}(Q)=\lim_{n\to\infty}4^{-n}
h_E(2^nQ).
\end{equation*}
There are explicit formul{\ae}{\ }for each of the
local heights (see \cite{silverman-1986}, and
\cite{silverman-1994}, \cite{everest-nyj-1999} for
an alternative approach).
For a prime $p$ where $Q$ has non-singular reduction,
\begin{equation}
\label{theneedleandthedamagedone}
\lambda_p(Q)=\textstyle\frac{1}{2}\log^{+}\vert x(Q)\vert_p.
\end{equation}
Notice in particular that if $x(Q)$ is integral at $p$
and $Q$ has non-singular reduction at $p$ then $\lambda_p(Q)=0$.
The singular reduction case is more involved, and to avoid
a major digression we deal only with {\sl split multiplicative
reduction} (see \cite[p. 362]{silverman-1994} for details on this):
the results all hold more generally but require passage to
extension fields.
In the split multiplicative case, the points on the
curve are isomorphic to the group ${\Bbb Q}_p^{*}/\ell^{\Bbb Z}$
where $\ell\in{\Bbb Q}_p^{*}$ has $\vert\ell\vert_p<1$. The
explicit formul{\ae}{\ }for
the $x$ and $y$ coordinates of a
non-identity point on the Tate curve
are given in terms of the parameter
$u\in{\Bbb Q}_p^{*}$ by
\begin{eqnarray*}
x&=&x_u=\sum_{n\in{\Bbb Z}}\frac{\ell^nu}{(1-\ell^nu)^2}-
2\sum_{n\ge1}\frac{n\ell^n}{(1-\ell^n)^2},\\
y&=&y_u=\sum_{n\in{\Bbb Z}}\frac{\ell^{2n}u^2}{(1-\ell^nu)^3}+
\sum_{n\ge1}\frac{n\ell^n}{(1-\ell^n)^2}.
\end{eqnarray*}
It is clear that $x_u=x_{u\ell}$ and $x_u=x_{u^{-1}}$.
Suppose $Q$ corresponds to the point $u\in{\Bbb Q}_p^{*}$
and assume, by invariance under multiplication by $\ell$,
that $u$ lies in the fundamental domain
$\{u\mid p^{-k}=\vert\ell\vert_p<\vert u\vert_p
\le 1\}$.
Then (by \cite{everest-nyj-1999} or
\cite{silverman-1994}),
\begin{equation*}
\lambda_p(Q)=\left\{
\begin{array}{ll}
-\log\vert 1-u\vert_p&\mbox{if }\vert u\vert_p=1,\\
-\frac{k}{2}\left(\frac{r}{k}-\left(\frac{r}{k}\right)^2\right)
&\mbox{if }\vert u\vert_p=p^{-r}<1.
\end{array}
\right.
\end{equation*}
Notice that for $\vert u\vert_p=1$,
the local height is non-negative, while if $\vert u\vert_p<1$
the local height is negative.
Also, these formula{\ae}{\ }extend to all of
$E(\Omega_p)$ by \cite{silverman-1994}
($\Omega_p$ is a fixed algebraic closure of ${\Bbb Q}_p$).
In \cite{dambros-thesis}, \cite{dambros-everest-miles-ward},
\cite{everest-ward-1998} and \cite{everest-ward-1999}, attempts
have been made to define dynamical systems whose topological entropy
is given by ${\hat{h}}(Q)$, the global canonical height of $Q$.
In the spirit of the algebraic case, and to reflect the fact that the
global canonical height is a sum of local canonical heights, one
looks to realize each local height as the entropy of a corresponding local
component. In \cite{dambros-thesis} and
\cite{dambros-everest-miles-ward} the elliptic adeles are
used. D'Ambros works over function fields and assumes that
the point $Q$ has everywhere non-singular reduction. In
\cite{dambros-everest-miles-ward} a similar non-singular reduction
assumption is made, together with an assumption that $Q$ lies
in a neighbourhood of the identity; there is also an
artificiality in the construction. Of particular interest is the
fact that the coincidence between periodic point
counts and division polynomials seen above holds
asymptotically
(cf. Remark \ref{broads}).
The extra freedom of sequential actions allows a different
approach to these problems, and gives a very clear
dynamically motivated description of the global canonical
height and the phenomenon of singular reduction. Now
the simple arithmetic structure of ${\Bbb P}^1({\Bbb Q})$
is replaced by the richer arithmetic of $E$.

\section{Duplication on elliptic curves}
\label{elliptics}

To fix notation, let $E$ be given in generalized Weierstrass form as
in (\ref{genweierstrassform}).
It follows from the shape of this equation that the denominator
of the $x$-coordinate of any rational point is a square.
Write $x(2^nQ)=\theta_n=a_n/b_n^2$, $b_n>0$ as a rational
in lowest terms.
\begin{remark}\label{bondingmaps}
It follows from the explicit formul{\ae}{\ }for
duplication that the
sequence of integers $(b_j)$ satisfies the strong divisibility property
$b_i|b_j\mbox{ for }i<j$,
that ensures the existence of well-defined transitional
maps $T_j^k:{\Bbb T}\to{\Bbb T}$ for $k\ge j$ with the property that
$T_k(x)=T_j^k(T_j(x))$ where $T_j(x)=b_jx$ mod $1$.
This brings the family of maps
$(T_j)$ closer to the iterates of a single transformation.
\end{remark}
According to Lemma \ref{noneoncompacts}, if we
are to realize the canonical height as the entropy
of a sequence of transformations on a compact space
(say the circle ${\Bbb T}$) then the rate must be
$r(n)=n$ essentially.
However, such systems cannot exhibit interesting entropies
unless they are of a very special shape.
\begin{lemma}\label{firstgresystemmodified} Let $X={\Bbb T}$,
and $T_j(x)=b_jx$ mod $1$ for $j>1$, where
$b_j\vert b_{j+1}$ for all $j\ge 1$, and
$b_{j+1}/b_{j}\to\infty$ as $j\to\infty$.
Then $h_{Bowen}^{r}({\bold T})=\infty$ for
$r(n)=n$ and is zero for
$r(n)/n\to\infty$.
\end{lemma}

\begin{pf} Let $B_{\epsilon}=(-\epsilon,\epsilon)$ and
think of the circle ${\Bbb T}$ as $[-1/2,1/2)$.
Notice that it is not possible to use the set
$B$ of Example \ref{simpleadelicexamples} since the group is compact.
For small $\epsilon$ (specifically, for
$\epsilon b_1<1$),
$B_{\epsilon}\cap T_1^{-1}B_{\epsilon}$
is a single interval, so
$$
\mu(B_{\epsilon}\cap T_1^{-1}B_{\epsilon})=
b_1^{-1}(2\epsilon).$$
However, for large $k$ the pre-image
under $T_k$ of $B_{\epsilon}$ meets
$B_{\epsilon}$ in a union of intervals:
$$
\mu(B_{\epsilon}\cap T_k^{-1}B_{\epsilon})=
\lceil b_k\epsilon\rceil b_k^{-1}(2\epsilon)+O(\epsilon/b_k).
$$
It follows that for large $k$,
$$
\mu\left(
\bigcap_{j=1}^{k}T_j^{-1}B_{\epsilon}\cap
T_k^{-1}B_{\epsilon}\right)=
\epsilon\mu\left(
\bigcap_{j=1}^{k}T_j^{-1}B_{\epsilon}\right)+O(\epsilon b_{k-1}/b_k)
$$
which shows the entropy is at least
$-\log(\epsilon)$. It follows that
$h_{Bowen}^{r}({\bold T})=\infty$ for $r(n)=n$.
The last case follows from Lemma \ref{noneoncompacts}.
\end{pf}
In the arithmetically simple case where the
$b_j$ are all powers of a single number, Lemma \ref{firstgresystem}
has a simpler proof. For example, if $b_j=2^{c_j}$, then
the factor map $x\mapsto\sum_{n=1}^{\infty}x_n2^{-n}$ from
the shift space $\Sigma=\{0,1\}^{\Bbb N}$ onto ${\Bbb T}$
intertwines doubling with the left shift.
The maps $T_j$ are then factors of powers of the shift:
$T_j=\sigma^{c_j}$, where $\sigma$ is the left shift
on $\Sigma$. If $B$ is now a cylinder set defined
on finitely many coordinates in $\Sigma$, then
it is clear that the sets
$T_j^{-1}B$ for distinct large $j$ are independent,
which shows that the topological entropy is infinite.

On the other hand, the non-compact analogue of this
system does exhibit interesting volume growth.
\begin{theorem}\label{firstgresystem}Let $r(n)=4^n$, $X={\Bbb R}$,
and $T_j(x)=b_jx$ for $j>1$ with the sequence $(b_n)$ 
defined by $x(2^nQ)=a_n/b_n^2$.
Then $h_{Bowen}^r({\bold T})=
\hat{h}(Q)$ for non-torsion $Q$.
\end{theorem}

\begin{pf} By a strong form of Siegel's theorem
(see \cite[p. 250]{silverman-1986}),
\begin{equation}\label{abssimilarinsize}
\lim_{n\to\infty}\frac{\log\vert a_n\vert}{2\log\vert b_n\vert}=1.
\end{equation}
Also,
\begin{equation}\label{defcanhgtsilv}
\lim_{n\to\infty}\frac{1}{r(n)}\log\frac{1}{2}
\max\{\vert a_n\vert,\vert b_n^2\vert\}=\hat{h}(Q)
\end{equation}
by \cite[Chap. VIII, Sect. 9]{silverman-1986}.
Thus $\vert b_n\vert\to\infty$ and
$\lim_{n\to\infty}\frac{1}{r(n)}\log b_n=\hat{h}(Q)$
by (\ref{abssimilarinsize}) and (\ref{defcanhgtsilv}).
It follows that
\begin{equation*}
\log\mu\left(
\bigcap_{j=1}^{n}T_j^{-1}B_{\epsilon}\right)=
-\log\max_{1\le j\le n}\{\vert b_j\vert\}+\log2\epsilon.
\end{equation*}
For any real sequence $(d_n)$ with
$
\frac{d(n)}{r(n)}\longrightarrow\omega\ge0,
$
\begin{equation}
\label{lateron}
\frac{\max_{1\le j\le n}\{d(j)\}}{r(n)}\longrightarrow\omega\ge0.
\end{equation}
It follows that $\frac{1}{r(n)}
\log\max_{1\le j\le n}\{\vert b_j\vert\}\to\hat{h}(Q)$
as required.
\end{pf}

\begin{theorem}\label{secondgresystem} Let $r(n)=4^n$, $X={\Bbb Q}_{\Bbb 
A}$,
and $T_n({\mathbf x})=\theta_n{\mathbf x}$
where $\theta_n={a_n}/{b_n^2}=x(2^nQ)$.
Then $h_{Bowen}^{r}({\bold T})=2\hat{h}(Q)$.
\end{theorem}

\begin{pf} It is enough to measure the volume growth of
the open set $B=(-1,1)\times\prod_{p<\infty}{\Bbb Z}_p.$
At the infinite place,
we need a bound on $\max_{1\le n\le N}\{\vert\theta_n\vert\}$,
and this is provided by elliptic transcendence theory
(see \cite{david-1991}).
The minimum distance of $nQ$ from the identity
on ${\Bbb C}/L$ is bounded below
by $n^{-A}$ for some $A=A(E,Q)>0$.
The size of the $x$-coordinate
is approximately the inverse square of this quantity. Since we are
running through the powers of $2$ only, this gives an upper bound for
$\max_{1\le n\le N}\{\vert\theta_n\vert\}$ of the shape $C^{N}$. Thus, if
\begin{equation*}
\bigcap_{j=1}^{N}T_j^{-1}B=B_{N,\infty}\times\prod_{p<\infty}
B_{N,p},
\end{equation*}
the measure of $B_{N,\infty}$ is $O(C^{N})$.
For the finite places, the sequence $(b_n)$ -- and hence
$(b_n^2)$ -- has a very strong divisibility property:
$b_i\vert b_{i+1}$ for all $i\ge1$ (by the
duplication formula).
Thus
\begin{eqnarray*}
\mu\left(B_{N,p}\right)&=&
\mu\left(\bigcap_{n=1}^{N}\left({a_n}/{b^2_{n}}\right)^{-1}
{\Bbb Z}_p\right)\\
&=&\min_{1\le n\le N}\left\{
\left\vert{a_n}/{b_n^2}\right\vert_p^{-1}\right\}\\
&=&\vert b_{N}\vert_p^2.
\end{eqnarray*}
It follows that
\begin{eqnarray*}
\log\mu\left(
\bigcap_{j=1}^{N}T_j^{-1}B\right)&=&
2\log\prod_{p<\infty}\vert b_{N}\vert_p
+O(\log C^N)\\
&=&-2\log\vert b_N\vert+O(N).
\end{eqnarray*}
So $h_{Bowen}^{r}({\bold T})=2\hat{h}(Q)$ as in the proof
of Theorem \ref{firstgresystem}.
\end{pf}

\section{A dynamical interpretation
of singular reduction}\label{dynintofbadred}

The systems described in Example \ref{solenoidsystemexamples}
have local entropies which sum to the global topological entropy.
Example \ref{simpleadelicexamples}
shows that the entropy of simple examples of sequences of transformations
on the adeles may not add up in an analogous way. In pursuit of
the connection between heights and entropy on elliptic curves,
a more substantial problem appears, preventing Theorems
\ref{firstgresystem} and \ref{secondgresystem} from decomposing
into local contributions. On the height side, it is still
the case that the global canonical height is a sum of local heights,
\begin{equation}\label{sumoflocalsinellipticcaseheight}
\hat{h}(Q)=\sum_{p\le\infty}\lambda_p(Q),
\end{equation}
(see \cite[App. C, Sect. 18]{silverman-1986}). When $p$ is a prime of
singular
reduction for the curve, or $p=\infty$, it is possible for the local
height $\lambda_p(Q)$ to be strictly negative. This means that
it certainly cannot represent the topological entropy of anything,
even in the sense of Definition \ref{mtseqentdef}.
In \cite{dambros-everest-miles-ward}, an approach to
interpreting the global height as the entropy of a dynamical
system is presented. Roughly speaking, since
(\ref{sumoflocalsinellipticcaseheight}) decomposes into
an expression for the global canonical height as the difference of
two non-negative quantities, it was suggested there that a
global system on the adeles might have a canonical factor, whose
quotient has the canonical height as entropy, and whose fibres carry the
other component of the entropy.

If $P=[x(P),y(P)]$ denotes a generic point on the curve $E$,
described by a generalized Weierstrass equation as before,
then
$x(nP)$ is a rational function of $x$ and $y$. In particular,
the denominator of that rational function is a polynomial
which vanishes on the $n$-torsion of $E$. This polynomial
can be used to generate a sequence of transformations
with more arithmetical subtlety.
Let $\psi_n$ denote the $n$th division polynomial of $E$
for $n\ge 1$ (see \cite[App. C]{everest-ward-1999},
\cite{silverman-1986}).
Thus, $\psi_n$ is an integral polynomial of degree $n^2-1$
and leading coefficient $n^2$ whose roots are exactly the
$x$-coordinates of all the non-identity points of order dividing
$n$ on $E$. It is well-known that $\psi_n(x)$ is always the
square of a polynomial in both $x$ and $y$ and, for odd $n$,
it is the square of a polynomial in $x$ alone (see
\cite[p. 105]{silverman-1986}).
Writing $q=a/b=x_Q$ for
the $x$-coordinate of a fixed rational point $Q$, define
\begin{equation*}
q_n=\vert b^{n^2-1}\psi_{n}(a/b)\vert\in{\Bbb Z}.
\end{equation*}
The remarks above show that $q_n$ is a square for all $n\ge1$.
Additionally, the sequence $(q_n)$ is a divisibility sequence
in the usual sense: $m\vert n$ implies $q_m\vert q_n$.
\begin{remark}
\label{broads}In the broad analogy being pursued,
the obvious candidate for the cardinality of periodic points is the sequence
$q_n$. However, if $f:X\to X$ is any bijection, then the
periodic points of $f$ must satisfy the combinatorial congruence
\begin{equation}\label{congruenceproperty}
0\le\sum_{d\vert n}\mu(n/d)\times\#\{x\in X\mid
f^d(x)=x\}\equiv 0\mbox{ mod }n, 
\end{equation}
for each $n\ge 1$.
Taking $E:y^2-y=x^3-x$ as the curve, and
$Q=(0,0)$ as the point, the sequence
$\vert\psi_n(0)\vert$ begins
$1,1,1,1,5,\dots$ which does not
satisfy (\ref{congruenceproperty}).
\end{remark}

These elliptic divisibility sequences $(q_n)$ were studied
in an abstract setting by Morgan Ward in a sequence
of papers - see \cite{mward-1948a} for the details.
Shipsey's thesis \cite{shipsey-thesis} contains more
recent applications of these sequences.

Define a sequence of non-negative
integers by $u_n^2=q_{2^n}$. If $Q$ is
not a torsion point then the terms of the sequence $(u_n)$ are
always non-zero. The divisibility of the sequence $(q_n)$
implies that
\begin{equation*}
u_1\vert u_2\vert u_3\vert\dots.
\end{equation*}
Define a sequence of transformations on ${\Bbb Q}_{\Bbb A}$
by
\begin{equation}\label{ontheadelenewsyst}
U_j({\bold x})=u_j^{-1}{\bold x}
\end{equation}
for $j\ge1$.
In Theorems \ref{firstgresystem} and \ref{secondgresystem}
the denominator of $\theta_n$ is responsible for
the volume growth, and hence the entropy.
These denominators may be thought of as evaluations
of the division polynomial (though in practice
a large amount of cancellation takes place).
Let $S$ denote the set of primes for which the point
$Q$ has singular reduction, and define the $S$-adeles to be
${\Bbb Q}_{S}=\prod_{p\in S}{\Bbb Q}_p$.
Write ${\bold U}_S$ for restriction of
$\bold U$ to ${\Bbb Q}_{S}$.
The local height of $Q$ is non-positive for each prime in
$S$, while for any prime $p$ dividing $b$, $Q$ has
non-singular reduction and the local height there is
$-\frac{1}{2}\log\vert 
b\vert_p$.

\begin{theorem}\label{goodbadreduction}
For the sequence of transformations {\rm(\ref{ontheadelenewsyst})}
and $r(n)=4^n$,

\noindent1. $h_{Bowen}^{r}
({\bold U})=\lambda_{\infty}(Q)+\frac{1}{2}\log\vert b\vert$,

\noindent2. $h_{Bowen}^{r}
({\bold U}_{{S}})=-\sum_{p\in S}\lambda_p(Q)\ge0$, and

\noindent3. $h_{Bowen}^{r}
(\bar{\bold U})=\hat{h}(Q)=\lambda_{\infty}(Q)+
\frac{1}{2}\log\vert b\vert+
\sum_{p\in S}\lambda_p(Q)$ where $\bar{\bold U}$ is the quotient
sequence of transformations
induced by
${\bold U}$ on ${\Bbb Q}_{\Bbb A}/{\Bbb Q}_{S}$.
\end{theorem}
Notice that the first formula is an analogue of Yuzvinskii's
formula. Theorem \ref{goodbadreduction} will be proved later.

\begin{corollary} If $Q$ has everywhere non-singular reduction
then
$$
h_{Bowen}^r({\bold U})=\hat{h}(Q).$$
If $Q$ has singular reduction at $p \in S$ then, with $\bar{\bold U}$
as before,
$$
h_{Bowen}^{r}(\bar{\bold U})=\hat h(Q).$$
\end{corollary}

Define $\epsilon_p(Q)$ to be $1$ if $\lambda_p(Q)\ge 0$
and $-1$ if $\lambda_p(Q)<0$. This map has the following
properties.

\noindent1. If $Q$ is integral, then $\epsilon_{\infty}(Q)=1$
(see comments after (\ref{firststepininfinitecase})).

\noindent2. The set of primes $p$ for which
$\epsilon_p(Q)=-1$ is finite.

\noindent3. There is a finite-index subgroup in $E({\Bbb Q})$
on which $\epsilon_p(Q)=1$ for all $p\in S$
(and therefore for all finite $p$) -- see
\cite[Sect. 6.2]{everest-ward-1999} or
\cite[Sect. 5]{dambros-everest-miles-ward}.

\noindent4. For all $Q$ in a neighbourhood of the
identity, $\epsilon_p(Q)=1$.

\begin{theorem}\label{flippingcaselocal}
For the sequence of transformations on ${\Bbb Q}_p$ defined by
$T_j(x)=q_j^{\epsilon_p(Q)}x$ for $j\ge1$,
where $Q\in E({\Bbb Q})$ is a non-torsion point, $q=x(Q)$,
$q_j^2=\vert\psi_{j}(q)\vert$
and $r(n)=n^2$,
$h_{Bowen}^{r}({\bold T})=\epsilon_p(Q)\lambda_p(Q)$.
\end{theorem}

\begin{pf} There are three cases to consider.
If $p=\infty$,
we claim firstly that
\begin{equation}
\label{firststepininfinitecase}
\lim_{N\to\infty}N^{-2}\log\vert\psi_N(q)\vert=2\lambda_{\infty}(Q).
\end{equation}
Notice that this explains the first of the properties
of $\epsilon_{\infty}$ above: if $Q$ is integral,
then the left-hand side of (\ref{firststepininfinitecase}) is
non-negative for all $N$.
Formula (\ref{firststepininfinitecase}) was proved in
\cite[Theorem 6.18]{everest-ward-1999}; the proof is
sketched here because it is similar to the singular reduction
case.
Take $G=E({\Bbb C})$ and consider the elliptic Jensen
formula
\begin{equation}\label{ejf}
\int_G\log\vert x(P)-x(Q)\vert d\mu_G(P)=2\lambda_{\infty}(Q)
\end{equation}
where $\mu_G$ is the normalized Haar
measure on $G$ (see \cite{everest-brid-1996}).
The points of $N$-torsion are dense
and uniformly distributed in $E({\Bbb C})$ as $N\to\infty$,
so the limit sum over the torsion points will tend to the
integral when the integrand is continuous.
The only potential problem arises from torsion points
close to $Q$: by \cite{david-1991}, for $x=x(P)$
with $NP=0$,
$\vert x-x(Q)\vert>N^{-C}$ for
some $C>0$ which depends on $E$ and $Q$
only.
This inequality is enough to imply
that the Riemann sum given by the $N$-torsion
points for $\log\vert x(P)-x(Q)\vert$ converges,
which gives (\ref{firststepininfinitecase}).
Now $q_n^2=\vert\psi_n(q)\vert$, so
\begin{equation*}
\log\mu\left(\bigcap_{j=1}^{N}T_j^{-1}B_{\epsilon}\right)
=-\log e_{N}+\log\epsilon,
\end{equation*}
where $e_N=\max_{1\le j\le N}\{q_n^{\epsilon_p(Q)}\}$,
so using (\ref{firststepininfinitecase}) gives
\begin{equation*}
\lim_{N\to\infty}N^{-2}\log e_N=\epsilon_{\infty}(Q)
\lambda_{\infty}(Q).
\end{equation*}

Assume that $p$ is a prime of singular reduction.
If $\vert x(Q)\vert_p>1$ then $Q$ has non-singular reduction
at $p$ and the result follows from the final case
below. Assume therefore that $\vert x\vert_p\le 1$,
and use the parametrisation of the curve
described in Section \ref{heightssection}.
The explicit formul{\ae}
of that section show that
the local height is non-positive.
The points of order dividing $N$ on the Tate curve
are precisely those of the form $\zeta^i\ell^{j/N}$,
$1\le i,j\le N$, where $\zeta\in\Omega_p$ denotes a fixed,
primitive $N$th root of unity in $\Omega_p$.
We claim that
\begin{equation}\label{andthereisnoplace}
\lim_{N\to\infty}N^{-2}\log\vert\psi_N(q)\vert_p=2\lambda_{p}(Q);
\end{equation}
this gives another proof that the local height is non-positive
at a point which is $p$-integral, where $p$ is a prime of singular
reduction.
Let $G$ denote the closure of the torsion points:
$G$ is not compact, so the $p$-adic elliptic Jensen
formula cannot be used. Instead we use a variant of
the Shnirelman integral: for $f:E(\Omega_p)\to{\Bbb R}$ define
the elliptic Shnirelman integral to be
\begin{equation*}
\int_Gf(Q)\mbox{d}Q=
\lim_{N\to\infty}N^{-2}
\sum_{N\tau=0}f(\tau)
\end{equation*}
whenever the limit exists.

We claim firstly that for any $P\in E({\Bbb Q}_p)$, the Shnirelman
integral
\begin{equation}\label{imgoingto}
\int_G\lambda_p(P+Q)\mbox{d}Q
=S(E)\mbox{ exists and is independent of $P$.}
\end{equation}
First assume that $P$ is the identity.
Using the explicit formula for the local height
gives
\begin{equation}\label{playasongforme}
-N^{-2}\sum_{i=1}^{N-1}\log\vert1-\zeta^i\vert_p
-N^{-2}\sum_{i=0}^{N-1}\sum_{j=1}^{N-1}
\frac{k}{2}\left(\frac{j}{N}-\left(\frac{j}{N}\right)^2\right).
\end{equation}
The first sum is bounded by $\log\vert N\vert_p/N$, which
vanishes in the limit; the second sum converges to
$-\frac{k}{12}$.

For the general case, let $P$ correspond to the point
$v$ on the multiplicative Tate curve.
If for some large $N$ no $j$ has $\vert\ell^{j/N}v\vert_p=1$
then the analogous sum to (\ref{playasongforme})
is close to $-\frac{k}{12}$ by the same argument.
Assume therefore that there is a $j$ with this property.
Then the first sum in (\ref{playasongforme}) is replaced
by
\begin{equation}\label{imnotsleepy}
-N^{-2}\sum_{i=0}^{N-1}\log\vert
1-\ell^{j/N}v\zeta^i\vert_p
-N^{-2}\log\vert1-(\ell^{r}v)^N\vert_p,
\end{equation}
where $r=j/N$ only depends on $v$.
By $p$-adic elliptic transcendence theory (see
\cite{david-1991}), there is
a lower bound for $\log\vert1-(\ell^rv)^N\vert_p$ of the
form $-(\log N)^A$, where $A$ depends on $E$ and $v=v(P)$
only. It follows that the first sum vanishes in
the limit as before. The second sum in
(\ref{playasongforme}) is simply rearranged under
rotation by $v$, so converges to $-\frac{k}{12}$ as before.
This proves (\ref{imgoingto}).

The claimed limit (\ref{andthereisnoplace})
now follows by taking the elliptic Shnirelman
integral of both sides of the parallelogram law
(\ref{paralawlocal}) and noting that equation (\ref{imgoingto})
shows that three terms cancel to leave the required limit.
Consider
$$
\log\mu\left(
\bigcap_{j=1}^{N}T_j^{-1}{\Bbb Z}_p\right)=
-\log f_N+\log\epsilon,
$$
where $f_N=\max_{1\le j\le N}\{\vert q_j\vert_p\}$.
Dividing by $N^2$ and taking the limit gives the result
as in the case $p=\infty$.

Finally, assume that $Q$ has non-singular reduction at $p$.
This is the easiest case. If $\vert x(Q)\vert_p=\vert q\vert_p>1$,
then
\begin{equation*}
\log\mu\left(
\bigcap_{j=1}^{N}T_j^{-1}B_{\epsilon}\right)=
-\log f_N+\log\epsilon,
\end{equation*}
where $f_N=\max_{1\le j\le N}\{\vert q_n\vert_p\}$.

If $(p,n)=1$ then $|\psi_n(q)|_p=|q|_p^{n^2-1}$.
If $(p,n)\neq 1$ then $(p,n-1)=1$. It follows
that
$$
|q|_p^{(N-1)^2-1}\leq f_N \leq |q|_p^{N^2-1}.
$$
Therefore
\begin{equation*}
\lim_{N\to\infty}N^{-2}\log\mu\left(
\bigcap_{j=1}^{N}T_j^{-1}B_{\epsilon}\right)=
\frac{1}{2}\log\vert q\vert_p=\lambda_p(Q)
\end{equation*}
by the explicit formula
(\ref{theneedleandthedamagedone}).
If $\vert q\vert_p\le 1$ then $q$ is a $p$-adic integer
and $\epsilon_p(Q)=1$. In this
case $\bigcap_{j=1}^{N}T_j^{-1}B_{\epsilon}=B_{\epsilon}$,
so there is no contribution to the entropy.
\end{pf}

\begin{pf} (of Theorem \ref{goodbadreduction})
For $p=\infty$
there can be no entropy contribution
for the sequence
$U_j({\bold x})=u_j^{-1}{\bold x}$, since $u_n$ is an integer sequence.
For $p$ finite, recall that $u_1\vert u_2\vert u_3\vert\dots$.
It follows that
\begin{equation*}
\bigcap_{j=1}^{N}U_j^{-1}\left(\prod_{p<\infty}{\Bbb Z}_p\right)
=u_{N}\left(\prod_{p<\infty}{\Bbb Z}_p\right),
\end{equation*}
which has measure $u_{N}^{-1}$.
The first result follows, since
$\log u_{N}=\frac{1}{2}(4^{N}-1)\log\vert b\vert+
\frac{1}{2}\log\vert\psi_{2^{N}}(q)\vert,$
by (\ref{firststepininfinitecase}).
The second follows at once by using ${\Bbb Q}_S$
and
(\ref{andthereisnoplace}).
For the
third part of the theorem, the calculation is the
same except that we are adrift by $\frac{1}{2}\sum_{p\in S,
(p,b)=1}\log|u_N|_p$. It follows from the proof of
Theorem \ref{flippingcaselocal}
that the entropy is adjusted by the contribution
of the local heights where $Q$ has singular
reduction.
\end{pf}

\section{Morphisms}
\label{goldstine}

Following the results of Call, Goldstine, Morton and Silverman
(see \cite{call-goldstine-1997}, \cite{morton-silverman-1994}),
we are able to place some of the results above in a broader
arithmetic-geometric context.
Recall that $F:{\Bbb P}^1\to{\Bbb P}^1$ is a morphism
of degree $d$ if it is given by $2$ homogeneous
polynomials of degree $d$, having only the point
$(0,\dots,0)$ as a common zero.
Write $[x,y]$ for the projective
coordinates in ${\Bbb P}^1$ and $z=x/y$. Then we
can think of a morphism $F$ as a single rational function $f$
in $z$. In this case, we will identify $F([z,1])$ with $[f(z),1]$
in the sequel.
For example, $f(z)=z^2$ is a morphism of degree $2$;
for $a,b\in{\Bbb Q}$ with $4a^3+27b^2\neq0$,
\begin{equation}\label{doublingonsphere}
f(z)=\frac{z^4-2az^2-8bz+a^2}{z^3+az+b}
\end{equation}
is a morphism of degree $4$.

The second example is the doubling map on
an elliptic curve.
Iteration
of rational maps has been
widely studied (see
\cite{beardon-1991} for an overview) and we
will show how our approach to dynamics can be developed in this
wider context. The rational map
(\ref{doublingonsphere}) on the Riemann
sphere already appears in both contexts:
in \cite[pp. 73-79]{beardon-1991}, Beardon gives this
map as the classic
example of a rational
function with empty Fatou set. The map appears in our work as
the map generating the sequential actions having entropies which
agree with classical notions of height.
Now recall the connection between heights
and complex dynamics. If
$F:{\Bbb P}^N\to{\Bbb P}^N$ is a morphism of
degree $d$ defined over ${\Bbb Q}$ then it
has an associated canonical height $\hat h_f$ with the properties

\noindent1. $\hat h_f(f(q)) = d\hat h_f(q)$ for any $q\in{\Bbb P}^N({\Bbb Q})$;

\noindent2. $q$ is pre-periodic if and only if $\hat h_f(q) = 0$.

A point $q$ is called {\sl pre-periodic} if the
orbit $\{f^n(q)\}$ is finite. In the elliptic example,
and working
over ${\Bbb C}$, the set of ${\Bbb C}$-pre-periodic points is
precisely the set of ${\Bbb C}$-torsion points. This accounts for
the example in \cite{beardon-1991} because the set of
${\Bbb C}$-torsion points
is dense in the complex elliptic curve.
Put differently,
the set of pre-periodic points in ${\Bbb P}^1({\Bbb C})$ under the rational
map is dense in ${\Bbb P}^1({\Bbb C})$.
Recent work in \cite{call-goldstine-1997}
and \cite{morton-silverman-1994} has decomposed
the global canonical height $\hat h_f$ into a sum of local
heights,
\begin{equation*}
\hat h_f(q)=\sum_p\lambda_{f,p}(q).
\end{equation*} 
These results, besides yielding beautiful
formul\ae{\ }have been used to give good bounds for the number
of pre-periodic points in certain cases.  

\begin{theorem}\label{localheightislocalentropy} Let $F$ denote
a morphism on ${\Bbb P}^1$ defined over $\Bbb Q$. Let 
$[q,1]$
denote a finite point of ${\Bbb P}^1(\Bbb Q)$.
Assume that the corresponding rational function
$f$ is actually a polynomial.
Then the iterates of $f$ on $q$ give a sequence
of rational numbers $f_n=f^n(q)$. The diagonal sequential
action on the adeles $T_n({\bold x})=
f_n{\bold x}$ has the following properties:

\noindent1. $h({\bold T})=\hat h_f(q)$

\noindent2. If ${\bold T}_p$ denotes the restriction of
${\bold T}$ to ${\Bbb Q}_p$
then $h({\bold T}_p)=\lambda_{f,p}(q)$.
\end{theorem}

\begin{pf}
Recall the following facts from
\cite{call-goldstine-1997}:
\begin{equation}\label{localheightisnicelimit}
\lambda_{f,p}(q)=\lim_{n\rightarrow \infty}\frac{1}{d^n}\lambda_p(f^n(q)),
\end{equation}
where $\lambda_p(a/b)$ is the local projective height 
$\lambda_p(a/b)=\log ^+|a/b|_p$.
The local height $\lambda_{f,p}(q)$ vanishes if and only if
$|f^n(q)|_p$ is bounded, for all $n$.
Finally, $q$ is pre-periodic if and only if
\begin{equation}\label{preperiodifflocalheightvanish}
 \lambda_{f,p}(q)=0
\text{ for all $p\le\infty$}.
\end{equation}
Note that properties
(\ref{localheightisnicelimit}) and 
(\ref{preperiodifflocalheightvanish}) only hold in the case
when $f$ is a polynomial: it is these properties 
which makes
the local heights so much easier to recognize as entropies. 

Write $f_n=a_n/b_n$, a rational in lowest
terms. Suppose
first that $|f^n(q)|_p$ is bounded
for some $p<\infty$. Then both the local
height and the local entropy are zero. If $|f^n(q)|_p$
is unbounded then
$\vert b_n\vert_p$ is eventually decreasing.
It follows that the volume growth
rate 
\begin{equation*}
\log \mu\left(\bigcap_{n=1}^{N}f_n^{-1}{\Bbb Z}_p\right)
=\log |b_{N}|_p
+O(1)=-\log^{+}\vert f_N\vert_p+O(1).
\end{equation*}
Dividing by $d^{N}$ and letting $N\rightarrow \infty$
shows this tends to the local canonical height, $\lambda_{f,p}(q)$.

For $p=\infty$,
\begin{equation*}
\log \mu\left(\bigcap_{n=1}^{N}f_n^{-1}B_{\epsilon}\right)
=-\max_{1\le n\le N}\log\vert f_{n}|
-\log2\epsilon,
\end{equation*}
which gives the result as in (\ref{lateron})
since $\log^{+}\vert f_n\vert/d^n\to
\lambda_{f,\infty}(q)\ge 0$.

The global case follows by combining the local ones:
$\bigcap_{n=1}^{N}f_n^{-1}B=B_{N,\infty}\times
\prod_{p<\infty}B_{N,p}$, and for all but finitely
many $p$, $B_{N,p}={\Bbb Z}_p$ for all $N$.
The volume growth rates of each of the (now finitely
many) local terms is the local height as shown above.
These sum to the global height.
\end{pf}

\begin{remark} Using a deeper form of Siegel's Theorem
due to Silverman \cite{silverman-1993},
it may be shown that the action $T_j({\bold x})=\theta_j
{\bold x}$ on ${\Bbb Q}_{\Bbb A}$, where
$\theta_j=f^j(q)$, has global entropy
equal to the global height for
any rational function $f$ fixing $\infty$.
\end{remark}

This result is much closer to the solenoid case we started
with. It also helps to put our elliptic results in a
better context. These show that the theorem is not true
in general when $f$ is a rational function; for example,
the local heights do not necessarily match up with the
local entropies. Nonetheless, there is a dynamical interpretation
of the values of the local heights. It would be of interest
to work out the general rational case along the lines of
the elliptic examples.

A deep uniformity in the behaviour of pre-periodic
points has been conjectured by Morton and Silverman.
\begin{conjecture} [Morton and Silverman] Let $F:{\Bbb P}^N({\Bbb C})\to
{\Bbb P}^N({\Bbb C})$ be a morphism of degree $d$ defined over
${\Bbb Q}$.
The number of pre-periodic points in ${\Bbb P}^1({K})$, where K
denotes an algebraic number field, is bounded by a constant
which depends on $N,d$ and $[K:Q]$ only.
\end{conjecture}

Proving this conjecture would have far-reaching
consequences: for example, it implies the
`Uniform Boundedness Conjecture' of
Mazur and Kamienny for torsion on elliptic curves
(and, more generally, on abelian varieties:
see \cite{merel-1996} for a proof of the former
and \cite{poonen-1997} for some discussion of the
latter).
In the examples we gave earlier, the `morphic' heights
correspond to the well known heights as follows.

\noindent1. For $f(z)=z^2$, $\hat h_f(q) = \log^{+}|q|$, the 
projective height. The local canonical height,
$\hat\lambda_{f,p}(q)=\log^{+}|q|_p,$
in agreement with the local component for the projective height.

\noindent2. For $f(z)=\frac{z^4-2az^2-8bz+a^2}{z^3+az+b}$,
the morphic height is precisely
the global canonical height and the local morphic heights
agree with the local canonical heights.

Notice now that the `circle' systems at the start can be
interpreted morphically. Given a rational $q$, the
sequence of squares generates a sequential action on the adeles
by defining $T_j({\bold x})=q^{2^j}{\bold x}$ (from repeated iteration of
$z\mapsto z^2$).
The sequential entropy agrees with the global morphic height,
and the local sequential entropies agree with the local morphic
heights.

\section{Heights, Periodic Points and the Julia set}\label{HMJS}

This paper has been about realizing elliptic 
heights -- and some morphic heights -- as entropies of
sequential
transformations, in
analogy with known circle results. The most convincing
elliptic examples rely upon recognizing local elliptic
heights as integrals. The space of integration turns out to
be the local curve, and this coincides
with the local Julia set of the associated
rational function. These examples give a good 
indication of how to
approach the more general problem of recognizing morphic heights
(both local and global) as entropies of sequences of 
transformations; namely, by recognizing the heights as
integrals over the Julia set.
We finish with an example and a theorem which illustrate this
connection between the Julia set and the 
morphic height. They suggest that Jensen's formula
(\ref{jensen}) could be
a fundamental stepping stone between heights and periodic
points in the morphic examples also.

\begin{theorem}\label{doughnut}
If $f(z)=az^d+\dots+a_0$ is a polynomial,
then for any $q\in{\Bbb C}$,
$$
\lambda_{f,\infty}(q)=\frac{1}{d-1}
\log\vert a\vert+\int_{J(f)}\log\vert x-q\vert\mbox{d}m(x),
$$
where $m$ is the maximal measure for $f$ on $J(f)$.
\end{theorem}

\begin{pf} Assume first that $q$ is not in the Julia set
of $f$.
The zeros of the
polynomial $f_n(x)=f^n(x)-x$ are precisely the solutions of the
equation $f^n(x)=x$.
Note that $d_n=\deg(f_n)=d^n$, where $d=\deg(f)$.
If $|f^n(q)|\to\infty$,
$\frac{1}{d_n}\log\vert f_n(q)\vert$ is approximately
$\frac{1}{d^n}\log\vert f^n(q)\vert$,
which converges to $\lambda_{f,\infty}(q)$ (the archimedean local height
of $q$ for the morphism associated to $f$).
If $|f^n(q)|$ is bounded, then the same is true
for $\vert f_n\vert$, so both expressions tend to
$\lambda_{f,\infty}(q)$.
Since $q$ lies in the open Fatou set,
$\log|x-q|$ is continuous on $J(f)$. Now
\begin{equation}
\label{iwasarrested}
\frac{1}{d_n}\log|f_n(q)|=\frac{1}{d_n}
\sum_{f^n(x)=x}\log|x-q|+\frac{1}{d_n}\log\vert B_n\vert,
\end{equation}
where the sum is over the $n$th `division points'
and
\begin{equation*}
B_n=a^{1+d+d^2+\dots +d^{(n-1)}}
\end{equation*}
is the
leading coefficient of $f^n(x)$.
Thus
\begin{equation}
\label{holdingpapers}
\frac{1}{d_n}\log\vert B_n\vert=
\frac{1}{d^n}\left(\frac{d^n-1}{d-1}\right)
\log\vert a\vert\to\frac{1}{d-1}\log\vert a\vert.
\end{equation}
Now it is known that
$$\frac{1}{d_n}\sum_{f^n(x)=x}\log|x-q|\to\int_{J(f)}
\log|x-q|\text{d}m(x),$$
where $m$ is the maximal invariant measure for $f$ restricted to the Julia
set (see \cite{lyubich-1982};
\cite{MR85m:58110b}).

It remains to show that the formula holds for $q\in J(f)$.
Without loss of generality, assume that $a=1$ (if
not, we may conjugate by a linear map to ensure this).
Since $J(f)$ has no interior, there is a sequence
$q_n\to q$ with $q_n\notin J(f)$.
Then $\log\vert x-q_n\vert\to
\log\vert x-q\vert$ for all $x\in J(f)\backslash\{q\}$.
Since $J(f)$ is bounded, $\log\vert x-q_n\vert$
and $\log\vert x-q\vert$ are uniformly bounded above
by $M$ say for $x\in J(f)\backslash\{q\}$.
So by Fatou's lemma
\begin{equation}
\label{bangbangmaxwell}
0=\lim_{n\to\infty}\int_{J(f)}\log\vert x-q_n\vert\mbox{d}m(x)
\le
\int_{J(f)}\log\vert x-q\vert\mbox{d}m(x)\le M.
\end{equation}
This shows that $x\mapsto\log\vert x-q\vert$ is
in $L^1(m)$.

Now $\vert x-f(q)\vert=\prod_{f(t)=x}\vert t-q\vert$,
so
\begin{eqnarray*}
\int_{J(f)}\log\vert x-f(q)\vert\mbox{d}m(x)&=&
\int_{J(f)}\sum_{f(t)=x}\log\vert t-q\vert\mbox{d}m(x)\\
&=&d\int_{J(f)}\log\vert x-q\vert\mbox{d}m(x)
\end{eqnarray*}
(the last equality follows from \cite{lyubich-1982}
or \cite[Theorem (d)]{MR85m:58110b}).
If $\int_{J(f)}\log\vert x-q\vert\mbox{d}m(x)>0$, then
the last equation contradicts (\ref{bangbangmaxwell}).
\end{pf}

\begin{example}\label{tchebysheffexample}
Consider
the Tchebycheff polynomial
of degree $d$,
$f(z)=T_d(z)=\cos(d\arccos(z))$.
The Julia set is the interval $J(f)=[-1,1]$.
The map $\phi:{\Bbb C}\to{\Bbb C}$ given by
$\phi(z)=\frac{1}{2}(z+z^{-1})$ is a semi-conjugacy from
$g:z\mapsto z^d$ onto $z\mapsto f(z)$, in other words, 
$f(\phi(z))=\phi(z^d)$.
Write $\psi$ for the branch of the inverse of
$\phi$ defined on $\{z\in{\Bbb C}\mid\vert z\vert>1\}$.
The canonical morphic height at
the infinite place is (for $q\notin J(f)$)
\begin{eqnarray*} 
\lambda_{f,\infty}(q)&=&\lim_{n\to\infty}\frac{1}{d^n}
\log^{+}\vert f^n(q)\vert\\
&=&\lim_{n\to\infty}\frac{1}{d^n}
\log^{+}\vert\phi g^n\psi(q)\vert\\
&=&\lim_{n\to\infty}\frac{1}{d^n}
\log^{+}\vert\frac{1}{2}\left(g^n\psi(q)+
\frac{1}{g^n\psi(q)}\right)\vert\\
&=&\lim_{n\to\infty}\max\left\{0,
\frac{1}{d^n}\log\vert g^n\psi(q)\vert\right\}\\
&=&\log^{+}\vert\psi(q)\vert.
\end{eqnarray*}
For $q\in J(f)$, the same formula holds since
$\hat{h}_{f,\infty}(q)=0$ there by \cite{call-goldstine-1997}
and $\log^{+}\vert\psi(q)\vert=0$ there by
a direct calculation.

Now by Jensen's formula, for any $q\in{\Bbb C}$,
\cite[Theorem 15.18]{rudin-real-and-complex-analysis},
\begin{eqnarray*}
\log^{+}\vert\psi(q)\vert&=&
\log2+\int_{{\Bbb S}^1}\vert
\phi(y)-q\vert\mbox{d}y\\
&=&
\log2+\int_{J(f)}\log\vert t-q\vert\mbox{d}m(t)
\end{eqnarray*}
since $m$ is the image under $\phi$ of the
maximal measure (Lebesgue) on the circle.
That is,
$$
\hat{h}_{f,\infty}(q)=\log2+\int_{J(f)}\log\vert t-q\vert\mbox{d}m(t).
$$

The constant $\log 2$ in $\hat{h}_{\infty}(q)$
may be explained in accordance with
Theorem \ref{doughnut}. The
leading coefficient of $T_d$ is $2^{d-1}$,
so $\frac{1}{d-1}\log\vert a\vert$ in this case is
exactly $\log2$.
\end{example}

A similar approach can be adopted in the case
of polynomials with connected Julia sets.
There the local conjugacy near $\infty$
extends to the whole domain of attraction of $\infty$,
which is the complement of the filled Julia set.

\bibliographystyle{amsplain}

\providecommand{\bysame}{\leavevmode\hbox to3em{\hrulefill}\thinspace}

\end{document}